\documentclass[11pt, final]{amsart}
\usepackage{amssymb}
\usepackage[mathscr]{eucal}

\theoremstyle{plain}
\newtheorem{theorem}{Theorem}

\newtheorem{corollary}[theorem]{Corollary}

\theoremstyle{remark}

\def\A{\mathscr{A}}
\def\E{\mathscr{E}}
\def\Re{\operatorname{Re}}
\def\Im{\operatorname{Im}}

\begin{document}
\title[]{On some automorphisms of the set of effects on Hilbert space}
\author{LAJOS MOLN\'AR}
\address{Institute of Mathematics\\
         Lajos Kossuth University\\
         4010 Debrecen, P.O.Box 12, Hungary}
\email{molnarl@math.klte.hu}
\thanks{  This research was supported from the following sources:\\
          1) Hungarian National Foundation for Scientific Research
          (OTKA), Grant No. T--030082 F--019322,\\
          2) A grant from the Ministry of Education, Hungary, Reg.
          No. FKFP 0304/1997}
\date{\today}
\begin{abstract}
The set of all effects on a Hilbert space has an affine structure (it is
a convex set) as well as a multiplicative structure (it
can be equipped with the so-called Jordan triple
product). In this paper we describe the corresponding automorphisms
of that set.
\end{abstract}
\maketitle

Let $H$ be a complex Hilbert space. Denote by $B(H)$ the
$C^*$-algebra of all bounded linear operators on $H$.
The operator interval $[0,I]$ of all positive operators in $B(H)$ which
are bounded by the identity $I$ is called the Hilbert space effect
algebra. This has important applications in quantum mechanics.
The effect algebra $[0,I]$ can be equipped with several algebraic
operations.
For example, one can define a partial addition on it. Namely, if $A,B
\in [0,I]$ and $A+B \in [0,I]$, then one can set $A\oplus B=A+B$.
This structure has been investigated in several papers (see \cite{CVLL1,
Gudder2, Gudder3} and the references therein).
Moreover, on $[0,I]$ there is a natural partial ordering $\leq $ which
comes from the usual ordering between the self-adjoint operators on $H$
and one can also define the operation of the so-called
orthocomplementation
by ${}^\perp: A \mapsto I-A$. The set of all effects on $H$ equipped
with
this ordering and orthocomplementation has been studied for example in
\cite{CVLL2}.
Next, $[0,I]$ is clearly
a convex subset of the linear space $B(H)$. So, one can consider the
operation of convex combinations on it. The set of all effects
with this structure has been investigated in \cite{Gudder3}, for
instance.
Finally, as for a mutliplicative operation on $[0,I]$, note that in
general $A,B \in [0,I]$ does not imply that $AB\in [0,I]$. However, we
all the time have $ABA\in [0,I]$. This multiplication which is a
nonassociative operation and sometimes called Jordan triple
product also appears in infinite dimensional holomorphy as well as in
connection with the geometrical properties of $C^*$-algebras.

Because of the importance of effect algebras,
it is a natural problem to study the isomorphisms of the mentioned
structures. The aim of this paper is to contribute to these
investigations.

The automorphisms of $[0,I]$ with the operation of partial addition were
described in \cite{CVLL1}. The automorphisms of the Hilbert space effect
algebra
equipped with the partial ordering $\leq $ and the orthocomplementation
${}^\perp$ were characterized in \cite{CVLL2} (see also \cite{Ludwig}).
The isomorphisms of $[0,I]$ as a convex subset of $B(H)$ were
investigated in \cite{Gudder3}. However, in that paper the authors
considered such affine functions (maps preserving convex combinantions)
which
are homogenous for the scalars in $[0,1]$. This means that they supposed
that
their affine bijections have the additional property that they send 0 to
0. As a corollary of our first theorem we describe the affine
isomorphisms of $[0,I]$ without this extra condition. In the second
theorem
we determine the automorphisms of $[0,I]$ equipped with the
Jordan triple product. It is worth mentioning that, as it turns
out from
our results, the linear and multiplicative structures of $[0,I]$ are
very closely related to each other.

Let us fix the notation and definitions that we shall use
throughout the paper. So, $B(H)$ and $B_s(H)$ denote
the $C^*$-algebra of all bounded linear operators on $H$ and the
$JB^*$-algebra of all bounded self-adjoint operators on $H$,
respectively. A self-adjoint idempotent $P$ in $B(H)$ is called a
projection.
A von Neumann algebra $\A$ on $H$ is said to be a factor if its center
is trivial, that is, it equals $\mathbb CI$ ($I$ is the identity on
$H$). Define the set $\E(\A)$ of all effects in $\A$ by
$\E(\A)=[0,I]\cap \A$.
If $\mathscr R_1, \mathscr R_2$ are *-algebras over the complex field,
then a linear map $\psi: \mathscr R_1 \to \mathscr R_2$ satisfying
$\psi(A)^*=\psi(A^*)$ $(A\in \mathscr R_1)$ is called
\begin{itemize}
\item[(i)] a Jordan *-homomorphism if $\psi(A)^2=\psi(A^2)$ $(A\in
\mathscr R_1)$;
\item[(ii)] a *-homomorphism if $\psi(A)\phi(B)=\psi(AB)$ $(A,B\in
\mathscr R_1)$;
\item[(iii)] a *-antihomomorphism if $\psi(A)\psi(B)=\psi(BA)$ $(A,B\in
\mathscr R_1)$.
\end{itemize}
If $X,Y$ are linear spaces over $\mathbb C$ and $C\subset X$ is a
convex
set, then the function $\psi: C \to Y$ is called affine if it satisfies
\[
\psi(\lambda x+(1-\lambda)y)=\lambda \psi(x) +(1-\lambda)\psi(y)
\]
for every $x,y\in C$ and $\lambda\in [0,1]$.

Our first result determines the affine automorphisms of $\E(\A)$ for any
factor $\A$.

\begin{theorem}
Let $\A$ be a factor. If $\phi: \E(\A) \to \E(\A)$ is a bijective affine
function, then there is an either *-automorphism or *-antiautomorphism
$\Phi$ of $\A$ such that
\[
\phi(A)=\Phi(A) \qquad (A\in \E(\A))
\]
or
\[
\phi(A)=\Phi(I-A) \qquad (A\in \E(\A)).
\]
\end{theorem}

\begin{proof}
First we recall the following fact whose proof requires only trivial
calculations.
Let $\phi: \E(\A) \to X$ be an affine function with $\phi(0)=0$ where
$X$ is a linear space. Define
\[
\Phi_1(A)=
\begin{cases}
0 & \text{ if }  A=0 \\
\| A\| \phi\bigl(\frac{A}{\| A\|}\bigr) & \text{ if } A\neq 0,
\enskip 0\leq A \in \A. \end{cases}
\]
Next let
\[
\Phi_2(A)=\Phi_1(A^+)-\Phi_1(A^-) \qquad (A^*=A\in \A),
\]
where $A^+$ and $A^-$ are the positive part and the negative part of
$A$, respectively. That is,
\[
A^+=(1/2)(|A|+A) \quad \text{ and } \quad
A^-=(1/2)(|A|-A).
\]
Finally, set
\[
\Phi(A)=\Phi_2(\Re A)+i \Phi_2(\Im A) \qquad (A\in \A),
\]
where $\Re A$ and $\Im A$ denote the real part and the imaginary part
of $A$, respectively.
Then $\Phi: \A \to X$ is the unique linear extension of $\phi$ from
$\E(\A)$ to $\A$.

Let $\phi :\E(\A) \to \E(\A)$ be an affine function. We assert that
$\phi$ is continuous in the norm topology. To see this, consider the
affine function
\[
\psi: A \longmapsto \phi(A)-\phi(0)
\]
on $\E(\A)$ which sends 0 to 0. Since its unique linear extension $\Psi:
\A \to \A$ has the property that $\Psi(A)+\phi(0) \in \E(\A)$ for every
$A\in \E(\A)$, we deduce that
\[
\| \Psi(A)\| \leq \| \Psi(A)+\phi(0)\| +\| \phi(0)\| \leq 2
\qquad (A\in \E(\A)).
\]
Clearly, every element $A$ of the unit ball of $\A$ can be written as
$A=A_1-A_2+i(A_3-A_4)$ for some $A_1,A_2,A_3,A_4\in \E(\A)$. It follows
that
$\Psi$ is bounded on the unit ball of $\A$. This implies that $\Psi$ and hence
$\phi$ are norm continuous.

Let now $\phi :\E(\A) \to \E(\A)$ be an affine bijection.
Then $\phi$ and it inverse are norm-continuous. Moreover, $\phi$
obviously preserves the extreme points of $\E(\A)$ which are exactly the
projections in $\A$.

We claim that $\phi(0)$ is either $0$ or $I$.
Let $P\neq 0,I$ be a projection in $\A$. If every projection in $\A$
commutes with $P$, then we obtain that every element of $\A$ commutes
with $P$ which, $\A$ being a factor, would imply that $P$ is a scalar
multiple of the identity but this is a contradiction. So, we can choose
a projection $Q$ in $\A$ which does not commute with $P$. Considering
the operator $U=I-2Q \in \A$ we get a unitary element in $\A$ which
does not
commute with $P$. So, we have $P\neq UPU^*$. In any von Neumann algebra
the unitary group is arcwise connected. Therefore, there is an
arc within the set of all projections in $\A$ connecting $P=IPI^*$ to
$UPU^*$. To sum up, every nontrivial projection in $\A$ can be connected
by an arc within the set of all projections to another projection
different from the first one. It is trivial that $0$ and $I$ can be
connected
only to themselves. Since $\phi$ is a homeomorphism of the set of all
projections in $\A$, we deduce that $\phi$ sends nontrivial projections to
nontrivial projections and hence we have either $\phi(0)=0$ or
$\phi(0)=I$. Clearly, we can assume without loss of generality that
$\phi(0)=0$ (otherwise, we consider the transformation $A\mapsto
I-\phi(A)$). Let $\Phi$ be the unique linear extension of $\phi$ onto
$\A$.
We already know that $\Phi$ is a bounded linear transformation which sends
projections to projections. It is a standard algebraic argument to
verify that $\Phi$ is a Jordan *-homomorphism
(see \cite[Remark 2.2]{BreSem} and use the spectral theorem of
self-adjoint operators together with the continuity of $\Phi$).

We assert that $\Phi$ is bijective.
If $\Phi(A)=0$, then we see that $\Phi_2(\Re A)=\Phi_2(\Im A)=0$.
Let $B,C$ denote the positive and negative parts of $\Re A$,
respectively, From $\Phi_2(\Re A)=0$ we infer that
$\Phi_1(B)=\Phi_1(C)$. Supposing that $B,C\neq 0$, this means that
\[
\| B\| \phi\biggl(\frac{B}{\| B\|}\biggr)=
\| C\| \phi\biggl(\frac{C}{\| C\|}\biggr).
\]
Using the homogenity of $\phi$ for the scalars in $[0,1]$, we conclude
that
\[
\phi\biggl(\frac{B}{\| B\|+\| C\|}\biggr)=
\phi\biggl(\frac{C}{\| B\|+\| C\|}\biggr).
\]
Since $\phi$ is injective, it follows that $B=C$ which gives us that
$\Re A =0$. Similarly, one can check that $\Im A=0$ is also true, so we
have $A=0$. Therefore, $\Phi$ is injective.
Since the range of $\Phi$ is a linear subspace of $\A$ which contains
$\E(\A)$
(recall that $\Phi$ is an extension of $\phi$), it follows that $\Phi$ is
surjective. So, $\Phi$ is a Jordan *-automorphism of $\A$.

It is well-known that every factor is a prime algebra. This means that
for any $A,B\in \A$, the equality
$A\A B=\{ 0\}$ implies that $A=0$ or $B=0$. Now, a
classical theorem of Herstein on Jordan homomorphisms \cite{Her} applies
to obtain that $\Phi$ is either a *-automorphism or a *-antiautomorphism
of $\A$. This completes the proof of the theorem.
\end{proof}

Taking into account the form of *-automorphisms and *-antiautomorphisms
of $B(H)$, we immediately have the following corollary.

\begin{corollary}
Let $\phi: [0,I] \to [0,I]$ be a bijective affine
function. Then there is an either unitary or antiunitary operator
$U$ on $H$ such that
\[
\phi(A)=UAU^* \qquad (A\in [0,I])
\]
or
\[
\phi(A)=U(I-A)U^* \qquad (A\in [0,I]).
\]
\end{corollary}

Our next result describes the automorphisms of $[0,I]$ equipped with the
Jordan triple product.

\begin{theorem}\label{T:triple1}
Suppose that $\dim H\geq 3$.
Let $\phi: [0,I] \to [0,I]$ be a bijective function satisfying
\begin{equation*}
\phi(ABA)=\phi(A)\phi(B)\phi(A) \qquad (A,B \in [0,I]).
\end{equation*}
Then $\phi$ is of the form
\begin{equation*}
\phi(A)=UAU^* \qquad (A\in [0,I]),
\end{equation*}
where $U$ is either a unitary or an antiunitary operator on $H$.
\end{theorem}

\begin{proof}
First observe that $\phi$ sends projections to projections. Indeed, if
$P\in B(H)$ is a projection, then we have $\phi(P)=\phi(P)^3$. Since
$\phi(P)$ is a positive operator, by the spectral mapping theorem we
obtain that $\sigma (\phi(P))\subset \{ 0,1\}$ and this proves that
$\phi(P)$ is a projection.

We next show that $\phi$ preserves the partial ordering $\leq$
among the projections. Let $P, Q\in B(H)$ be projections and suppose
that $P\leq Q$. Then we have $PQP=P$ which yields
$\phi(P)=\phi(P)\phi(Q)\phi(P)$. This implies that $\phi(P)\phi(Q)$ is
an idempotent. On the other hand, since $\phi(P)$ and $\phi(Q)$ are
projections, the norm of their product is not greater than 1. So,
$\phi(P)\phi(Q)$ is a contractive idempotent. It is well-known that this
implies that $\phi(P)\phi(Q)$ is a projection and hence, due to the
self-adjointness, it follows that $\phi(P)$ and $\phi(Q)$ are commuting.
Hence, we can compute
\[
\phi(P)=\phi(P)\phi(Q)\phi(P)=
\phi(Q)\phi(P)\phi(P)=
\phi(Q)\phi(P)
\]
which yields that $\phi(P) \leq \phi(Q)$.
Since $\phi^{-1}$ has the same properties as $\phi$, it follows that
$\phi$ preserves the ordering $\leq$ in both directions. In particular,
we obtain that $\phi(0)=0$, $\phi(I)=I$ and that $\phi$ preserves the
nonzero minimal projections, that is, the rank-one projections on $H$.

We claim that $\phi$ preserves also the orthocomplementation on the set
of
projections. To see this, we first show that $\phi$ preserves the mutual
orthogonality. Let $P, Q\in B(H)$ be projections such that $PQ=0$. Then
we have $0=\phi(PQP)=\phi(P)\phi(Q)\phi(P)$ which implies that
\[
0=\phi(P)\phi(Q)\phi(Q)\phi(P)=\phi(P)\phi(Q) (\phi(P)\phi(Q))^*.
\]
This gives us that $\phi(P)\phi(Q)=0$.
It follows that $\phi(P)+\phi(I-P)$ is a projection, say $\phi(Q)$.
Since $\phi(P), \phi(I-P)\leq \phi(Q)$ and $\phi$ preserves the ordering
in both directions, we infer that $P, I-P\leq Q$. This gives us that
$Q=I$ and, hence, $\phi(P)+\phi(I-P)=I$.
Therefore, $\phi$ preserves the orthocomplementation on the set of all
projections.
The form of such transformations, that is, the form of all bijections
of the set of all projections on a Hilbert space with dimension
not less than 3 which preserve
the order in both directions and the orthocomplementation, is well-known
(see, for example, \cite{CVLL1}).
Namely, there is an either unitary or antiunitary operator $U$ on $H$
such that
\[
\phi(P)=UPU^*
\]
for all projections $P$ on $H$.

We next prove that $\phi(\lambda P)=\lambda \phi(P)$ for every $\lambda
\in [0,1]$ and every rank-one projection $P$. In fact, in that case
we can compute
\[
\phi(\lambda P)=\phi(P(\lambda P)P)=\phi(P)\phi(\lambda
P)\phi(P)=f_P(\lambda )\phi(P)
\]
for some scalar $f_P(\lambda) \in [0,1]$ which follows from the fact
that $\phi(P)$ is of rank one.
We assert that $f_P$ is a multiplicative function.
If $\mu \in [0,1]$, then we have
\[
f_P(\lambda^2 \mu)\phi(P)=
\phi(\lambda^2 \mu P)=
\phi((\lambda P)(\mu P)(\lambda P))=
\]
\[
\phi(\lambda P)\phi(\mu P)\phi(\lambda P)=
f_P(\lambda)^2f_P(\mu)\phi(P)
\]
which implies that $f_P(\lambda^2 \mu)=f_P(\lambda)^2f_P(\mu)$.
Choosing $\mu=1$, it follows that $f_P(\lambda^2)=f_P(\lambda)^2$.
We next obtain that
$f_P(\lambda^2 \mu)=f_P(\lambda^2)f_P(\mu)$.
Since this holds for every $\lambda ,\mu \in [0,1]$, we conclude that
$f_P$ is multiplicative.
We now claim that $f_P$ does not depend on the rank-one projection $P$.
If $P,Q$ are rank-one projections which are not mutually orthogonal,
then $PQP\neq 0$ and we have
\[
f_Q(\lambda^2)\phi(PQP)=
f_Q(\lambda^2)\phi(P)\phi(Q)\phi(P)=
\phi(P)\phi(\lambda^2 Q)\phi(P)=
\]
\[
\phi(P(\lambda^2 Q)P)=
\phi((\lambda P) Q(\lambda P))=
\]
\[
\phi(\lambda P)\phi(Q)\phi(\lambda P)=
f_P(\lambda^2)\phi(PQP).
\]
This gives us that $f_P=f_Q$. If $P,Q$ are mutually orthogonal, then
there is a rank-one projection $R$ such that $PRP\neq 0$ and $RQR\neq
0$. Thus we have $f_P=f_R=f_Q$.
So, there is a multiplicative function $f:[0,1] \to [0,1]$ such that
\[
\phi(\lambda P)=f(\lambda )\phi(P)
\]
for every $\lambda \in [0,1]$ and rank-one projection $P$ on $H$.
We show that $f$ is also additive on $[0,1]$. To see this, for any unit
vector $x\in H$ denote by $P_x$ the rank-one projection onto the linear
subspace of $H$ spanned by $x$. Let $x,y\in H$ be mutually orthogonal
unit vectors and $\lambda,\mu \in [0,1]$ such that $\lambda^2+\mu^2=1$.
Then $z=\lambda x +\mu y$ is a unit vector. We compute
$\phi(P_ z(P_x+P_y)P_z)$ in two different ways. On the one hand,
since $P_ z(P_x+P_y)P_z=P_z$, we have $\phi(P_
z(P_x+P_y)P_z)=\phi(P_z)$.
On the other hand, we compute
\[
\phi(P_ z(P_x+P_y)P_z)=
\phi(P_ z)\phi(P_x+P_y)\phi(P_z)=
\phi(P_ z)(\phi(P_x)+\phi(P_y))\phi(P_z)=
\]
\[
\phi(P_ z)\phi(P_x)\phi(P_z)+\phi(P_ z)\phi(P_y)\phi(P_z)=
\phi(P_ z P_x P_z)+\phi(P_ z P_y P_z)=
\]
\[
\phi(\lambda^2 P_z)+\phi(\mu^2 P_ z)=
(f(\lambda^2) +f(\mu^2)) \phi(P_ z)
\]
where we have used the fact that $\phi$ is orthoadditive on the set of
all projections (this follows from the form of $\phi$ on that set).
Therefore, we have $f(\lambda^2)+f(\mu^2)=1=f(\lambda^2
+\mu^2)$. By multiplicativity, we obtain the additivity of $f$.
We claim that $f$ is in fact the identity on $[0,1]$.
Since $f$ maps into $[0,1]$, one can easily check that $f$ is
monotone increasing. Moreover, as $f(1)=1$, the additivity of $f$
implies
that $f(r)=r$ for every rational number $r$ in $[0,1]$. If $\lambda \in
]0,1[$
is arbitrary, then approximating $\lambda$ by rationals $r,s$ from below
and above, respectively, by the monotonity we can infer that
$f(\lambda)=\lambda$.

We already know the form of $\phi$ on the set of all projections.
It is easy to see that without loss of generality we can assume that
$\phi(P)=P$ holds for every projection $P$ and we then have to prove
that $\phi$ is the identity on the whole interval $[0,I]$.
But this is now easy. Indeed, let $A\in [0,I]$. Pick an arbitary
rank-one projection $P=P_x$, where $x\in H$ is a unit vector.
Then we compute
\[
P\phi(A)P=
\phi(PAP)=
\phi(\langle Ax,x\rangle P)=
\langle Ax,x\rangle \phi(P)=
\langle Ax,x\rangle P=
PAP.
\]
Since $P$ was arbitrary, we obtain $\phi(A)=A$ for every $A\in [0,I]$.
This completes the proof of the theorem.
\end{proof}

Since the Jordan algebra $B_s(H)$ of all self-adjoint operators
also plays very important role in the mathematical foundations of
quantum
mechanics, we were tempted to determine the automorphisms of the set
$B_s(H)$
equipped with the Jordan triple product. Observe that the following
theorem has the interesting corollary that every such automorphism is
automatically linear, so one can say that the linear structure
of $B_s(H)$ is completely determined by its multiplicative Jordan triple
structure.
We remark that the question when a multiplicative function is
necessarily
additive was investigated for associative rings (recall that our
structure is highly nonassociative) in the purely algebraic setting (see
\cite{Martindale}).

\begin{theorem}
Suppose that $\dim H\geq 3$.
Let $\phi: B_s(H) \to B_s(H)$ be a bijective function
(linearity is not assumed) satisfying
\begin{equation*}
\phi(ABA)=\phi(A)\phi(B)\phi(A) \qquad (A,B \in B_s(H)).
\end{equation*}
Then there is an either unitary or antiunitary operator $U$
on $H$ such that either
\begin{equation*}
\phi(A)=UAU^* \qquad (A\in B_s(H))
\end{equation*}
or
\begin{equation*}
\phi(A)=-UAU^* \qquad (A\in B_s(H)).
\end{equation*}
\end{theorem}

\begin{proof}
We first prove that $\phi(I)$ is either $I$ or $-I$. Since
\[
\phi(I)\phi(A)\phi(I)=\phi(A)
\]
for every $A\in B_s(H)$ and $\phi$ is surjective, it follows that
$\phi(I)^2=I$. Therefore, we have
\[
\phi(I)\phi(A)=
\phi(I)\phi(A)\phi(I)\phi(I)=
\phi(A)\phi(I).
\]
Since this holds for every $A\in B_s(H)$, by the surjectivity of $\phi$,
it follows that $\phi(I)$ is in the center of $B(H)$ and, consequently,
$\phi(I)$ is a scalar. This yields that either $\phi(I)=I$ or
$\phi(I)=-I$. Clearly,  the function $-\phi$ is a bijective mapping of
$B_s(H)$
satisfying the equation appearing in the statement. So, without loss
of generality we can assume that $\phi(I)=I$.

We prove that $\phi$ sends projections to projections. If $P$ is a
projection, then $\phi(P)$ is self-adjoint and we have
$\phi(P)^2=\phi(P)\phi(I)\phi(P)=\phi(PIP)=\phi(P)$ which shows that
$\phi(P)$ is an idempotent.

Now, we can follow the argument in the proof of
Theorem~\ref{T:triple1}. One can verify that
$\phi$ preserves the partial ordering $\leq $ in both directions and the
orthocomplementation on the set of all
projections. So, we have an either
unitary or antiunitary operator $U$ on $H$ such that
\[
\phi(P)=UPU^*
\]
for every projection $P$ on $H$. One can check that for every rank-one
projection $P$ there exists a function $f_P: \mathbb R \to \mathbb R$
such that $\phi(\lambda P)=f_P( \lambda)\phi(P)$ $(\lambda \in \mathbb
R)$. We next obtain that
$f_P(\lambda ^2\mu )=f_P(\lambda)^2 f_P(\mu)$ and
choosing $\mu=1$, this gives us that
$f_P(\lambda ^2)=f_P(\lambda)^2$. In particular, by the injectivity of
$f_P$, from
\[
f_P(-\lambda)^2=f_P((-\lambda)^2)=f_P(\lambda^2)=f_P(\lambda)^2
\]
we deduce that $f_P(-\lambda)=-f_P(\lambda)$. Since
$f_P(\lambda ^2\mu )=f_P(\lambda^2) f_P(\mu)$, we get that
$f_P$ is multiplicative.
One can next show that $f_P$ does not depend on the rank-one projection
$P$.
So, there is a multiplicative function $f:\mathbb R \to \mathbb R$ such
that
\[
\phi(\lambda P)=f(\lambda )\phi(P)
\]
for every real number $\lambda$ and rank-one projection $P$ on $H$.
As for the additivity of $f$, just as in the proof of our previous
theorem we get that
\[
f(t)=f(t\lambda^2)+f(t\mu^2)
\]
for every real $t$, where $\lambda^2+\mu^2=1$. To show that $f$ is
additive, it is enough to verify that $f(1)=f(t)+f(1-t)$ for every
real $t$. If $t\in [0,1]$, then we already know this. If $t\notin
[0,1]$, say $t<0$, then we can refer to the equality
\[
f(1-t)=f\biggl((1-t)\frac{1}{1-t}\biggr)+f\biggl((1-t)\frac{-t}{1-t}\biggr)
\]
what is known to be valid since
the numbers $\frac{1}{1-t}, \frac{-t}{1-t}$ belong to $[0,1]$ and their
sum is 1. So, we have
\[
f(1-t)=f(1)+f(-t)=f(1)-f(t).
\]
Consequently, we obtain that $f: \mathbb R \to \mathbb R$ is an
injective multiplicative and additive function.
This means that $f$ is a nontrivial ring endomorphism of $\mathbb R$.
It is well-known that $f$ is necessarily
the identity (anyway, this can be proved quite similarly to
the corresponding part of the proof of Theorem~\ref{T:triple1}).
Finally, one can complete the
proof of the statement just as in the case of our previous theorem.
\end{proof}

\end{document}